\documentclass[11pt]{amsart}

\usepackage{amssymb}
\usepackage{amsmath}
\usepackage{amscd}

\theoremstyle{plain}
\newtheorem{prop}{Proposition}[section]
\newtheorem{proposition}[prop]{Proposition}
\newtheorem{conjecture}[prop]{Conjecture}
\newtheorem{lemma}[prop]{Lemma}
\newtheorem{theorem}[prop]{Theorem}

\theoremstyle{definition}

\theoremstyle{remark}

\newtheorem{remark}[prop]{Remark}
\newcommand{\con}{\equiv}
\DeclareMathOperator{\End}{End}
\DeclareMathOperator{\Gal}{Gal}
\DeclareMathOperator{\Spec}{Spec}

\DeclareMathOperator{\Ker}{ker}
\DeclareMathOperator{\Tate}{Tate}
\DeclareMathOperator{\tor}{tor}

\newcommand{\vphi}{\varphi}
\newcommand{\isom}{\cong}
\newcommand{\tensor}{\otimes}

\newcommand{\Qbar}{\overline{\mathbf Q}}
\newcommand{\Q}{{\mathbf Q}}
\newcommand{\C}{{\mathbf C}}
\newcommand{\F}{{\mathbf F}}
\newcommand{\Z}{{\mathbf Z}}
\newcommand{\T}{{\mathbf T}}

\newcommand{\m}{{\mathfrak m}}
\newcommand{\ra}{\rightarrow}

\newcommand{\intersect}{\cap}
\newcommand{\hra}{\hookrightarrow}
\newcommand{\Jtwo}{\tilde{J}^{(2)}}
\DeclareFontEncoding{OT2}{}{} 
  \newcommand{\textcyr}[1]{%
    {\fontencoding{OT2}\fontfamily{wncyr}\fontseries{m}\fontshape{n}%
     \selectfont #1}}
\newcommand{\Sha}{{\mbox{\textcyr{Sh}}}}
\numberwithin{equation}{section}

\title{Modular Parametrizations of Neumann--Setzer Elliptic Curves}
\author{William Stein\and Mark Watkins}
\thanks{International Mathematical Research Notices 2004 \#27, 1395--1405}
\begin{document}
\begin{abstract}
Suppose~$p$ is a prime of the form $u^2+64$ for some integer~$u$,
which we take to be $3$~mod~$4$.  Then there are two Neumann--Setzer
elliptic curves~$E_0$ and $E_1$ of prime conductor~$p$, and both have
Mordell--Weil group $\Z/2\Z$.  There is a surjective map
$X_0(p)\xrightarrow{\pi} E_0$ that does not factor through any other
elliptic curve (i.e., $\pi$ is optimal), where $X_0(p)$ is the modular
curve of level~$p$.  Our main result is that the degree of $\pi$ is
odd if and only if $u \con 3\pmod{8}$.  We also prove the
prime-conductor case of a conjecture of Glenn Stevens, namely that
that if~$E$ is an elliptic curve of prime conductor~$p$ then the
optimal quotient of $X_1(p)$ in the isogeny class of~$E$ is the curve
with minimal Faltings height.  Finally we discuss some conjectures and
data about modular degrees and orders of Shafarevich--Tate groups of
Neumann--Setzer curves.
\end{abstract}
\maketitle
\section{Introduction}\label{sec:intro}
Let~$p$ be a prime of the form $u^2+64$ for some integer~$u$, which we
take to be~$3$ modulo~$4$.  Neumann and Setzer \cite{neumann, setzer}
considered the following two elliptic curves of conductor~$p$ (note
that Setzer chose $u\con 1\pmod{4}$ instead):
\begin{eqnarray}
 E_0:& y^2+xy &=x^3-\frac{u+1}{4}x^2+4x-u,\label{eqn:e0}\\
 E_1:& y^2 + xy &= x^3 -\frac{u+1}{4}x^2 -x.\label{eqn:e1}
\end{eqnarray}
For $E_1$ we have $c_4=p-16$ and $c_6=u(p+8)$ with $\Delta=p=u^2+64$,
while for $E_0$ we have $c_4=p-256$ and $c_6=u(p+512)$ with
$\Delta=-p^2$.  Thus each $E_i$ is isomorphic to a curve of the form
$y^2 = x^3 - 27c_4x - 54c_6$ for the indicated values of $c_4$ and
$c_6$.  The curves $E_0$ and $E_1$ are $2$-isogenous and one can show
using Lutz-Nagell and descent via $2$-isogeny that
$$E_0(\Q)=E_1(\Q)=\Z/2\Z.$$  Moreover, if~$E$ is
{\em any} elliptic curve over~$\Q$ of prime conductor with a rational
point of order~$2$ then~$E$ is a Neumann--Setzer curve or has
conductor~$17$ (see \cite{setzer}).

Let $X_0(p)$ be the modular curve of level~$p$.
By \cite{wiles} there is a surjective map
$\pi: X_0(p)\to E_0$, and
by \cite[\S5, Lem.~3]{mestre-oesterle} we may choose~$\pi$ to
be optimal, in the sense that~$\pi$
does not factor through any other elliptic curve. 
The {\em modular degree} of $E_0$ is $\deg(\pi)$.

We prove in Section~\ref{sec:thm} that the modular degree of~$E_0$ is
odd if and only if $u\con 3\pmod{8}$.  Our proof relies mostly on
results from \cite{mazur}.  In Section~\ref{sec:stevens} we show that
$E_1$ is the curve of minimal Faltings height in the isogeny class
$\{E_0,E_1\}$ of $E_1$ and prove that $E_1$ is an optimal quotient
of~$X_1(p)$, which is enough to prove the prime-conductor case of a
conjecture of \cite{stevens} (this case is not covered by the
results of \cite{vatsal}).  Finally, in Section~\ref{sec:conj}
we give evidence for our conjecture that there are infinitely many
elliptic curves with odd modular degree, and give a conjectural
refinement of Theorem~\ref{thm:main}.  We also present some data about
$p$-divisibility of conjectural orders of Shafarevich--Tate groups of
Neumann--Setzer curves.

\subsection{Notation}\label{sec:notation}
Let~$p$ be a prime and~$n$ be the
numerator of $(p-1)/12$.  

We use standard notation for modular forms, modular curves, and Hecke
algebras, as in \cite{diamond-im} and \cite{mazur}.  In particular, let
$X_0(p)$ be the compactified coarse moduli space of elliptic curves
with a cyclic subgroup of order~$p$. Then $X_0(p)$ is an algebraic
curve defined over~$\Q$.  Let $J=J_0(p)$ be the Jacobian of $X_0(p)$,
and let $\T=\Z[T_2, T_3,\ldots]\subset\End(J)$ be the Hecke algebra.
Also, let $X_1(p)$ be the modular curve the classifies isomorphism
classes of pairs $(E,P)$, where $P\in E$ is a point of order~$p$.

To each newform $f\in S_2(\Gamma_0(p))$, there is an associated
abelian subvariety $A=A_f\subset J_0(p)$.  We call the kernel $\Psi_A$
of the natural map $A\hra J \ra A^{\vee}$ the {\em modular kernel}.
For example, when~$A$ is an elliptic curve, this map is induced by pullback
followed by push forward on divisors and $\Psi_A$ is multiplication by
$\deg (X_0(p) \to A)$.  The {\em modular degree} of $A$ is the
square root of the degree of $\Psi_A$.  This definition makes
sense even when $\dim(A)>1$, since the degree of a polarization
is the square of its Euler characteristic, hence a
perfect square (see \cite[\S{}16, pg.~150]{mumford}).
If $I\subset\T$ is an ideal, let 
$$A[I] = \{ x \in A(\Qbar) : Ix = 0\}
\qquad\text{and}\qquad
A[I^\infty] = \bigcup_{n>0} A[I^n].$$

\subsection{Acknowledgements}
The authors would like to thank the American Institute of Mathematics
 for hospitality while they worked on this paper, the National Science
 Foundation for financial support, and Matt Baker, Frank Calegari, and
 Barry Mazur for helpful conversations.

\section{Determination of the Parity of the Modular Degree}\label{sec:thm}
Let~$p$, $E_0$, $J$ and $n$ be as in Section~\ref{sec:intro}, and
fix notation as in Section~\ref{sec:notation}.
In this section we prove the following theorem.

\begin{theorem}\label{thm:main}
The modular degree of~$E_0$ is odd if and only if $u \con 3\pmod{8}$.
\end{theorem}

In order to prove the theorem we deduce seven lemmas using techniques
and results from \cite{mazur}.

Let~$m$ be the modular degree of $E_0$, and let 
$$
  B=\ker(J \xrightarrow{\pi} E_0).
$$
\begin{lemma}
We have $m^2 = \#(B\intersect E_0)$.
\end{lemma}
\begin{proof}
  As mentioned in Section~\ref{sec:notation}, the composition $E_0\to
  J \to E_0$ is multiplication by the degree of $X_0(p)\to E_0$, i.e.,
  multiplication by the modular degree of~$E_0$.  The lemma follows
  since multiplication by~$m$ on $E_0$ has degree~$m^2$.
\end{proof}

The {\em Eisenstein ideal} $\mathcal{I}$ of~$\T$ is the ideal
generated by $T_\ell-(\ell+1)$ for $\ell\neq p$ and $T_p -1$.  By
hypothesis, there is a Neumann--Setzer curve of conductor~$p$, which
implies that the numerator~$n$ of $(p-1)/12$ is even (we do the
elementary verification that this numerator is even in the proof of
Theorem~\ref{thm:main} below).  As discussed in
\cite[Prop.~II.9]{mazur}, the {\em $2$-Eisenstein
prime}~$\m=(2)+\mathcal{I}$ of~$\T$ is a maximal ideal of~$\T$, with
$\T/\m \isom \Z/2\Z$.

\begin{lemma}\label{lem:mtwo}
We have
  $E_0[\m] = E_0[2]$.
\end{lemma}
\begin{proof}
  By \cite[Prop.~II.11.1, Thm.~III.1.2]{mazur}, the Eisenstein ideal
  $\mathcal{I}$ annihilates $J(\Q)_{\tor}$, so~$\m$ annihilates
  $J(\Q)_{\tor}[2]$.  Since $J(\Q)_{\tor}$ is cyclic of order~$n$ (by
  \cite[Thm.~III.1.2]{mazur}), $J(\Q)_{\tor}[2]$ has order~$2$, so
  $J(\Q)_{\tor}[2]= E_0(\Q)_{\tor}[2]$, hence $E_0(\Q)[\m]\neq 0$.
  The Hecke algebra $\T$ acts on $E_0$ through $\End(E_0)\isom \Z$, so
  each element of $\T$ acts on $E_0$ as an integer; in particular, the
  elements of~$\m$ all act as multiples of~$2$ (since $E_0[\m]\neq 0$
and $2\in\m$), so $E_0[\m]=E_0[2]$ since $2\in\m$.
\end{proof}

\begin{lemma}\label{lem:powkill}
  Suppose $A\subset J_0(p)$ is a $\T$-stable abelian subvariety and
  $\wp\subset\T$ is a maximal ideal such that $A[\wp^{\infty}]\neq 0$.
  Then $A[\wp]\neq 0$.  Also $A[\wp^\infty]$ is infinite.
\end{lemma}
\begin{proof}
  Arguing as in \cite[\S{}II.14, pg.~112]{mazur}, we see that for
  any~$r$, $A[\wp^r]/A[\wp^{r+1}]$ is isomorphic to a direct sum of
  copies of $A[\wp]$.  If $A[\wp]=0$, then since $A[\wp^\infty]\neq
  0$, there must exist an~$r$ such that $A[\wp^r]/A[\wp^{r+1}]\neq 0$.
  But $A[\wp^r]/A[\wp^{r+1}]$ is contained in a direct sum of
  copies of~$A[\wp]=0$, which is a contradiction.
    
    To see that $A[\wp^\infty]$ is infinite, note that if $\ell$ is
    the residue characteristic of~$\wp$ and $\Tate_\ell(A)$ is the
    Tate module of~$A$ at $\ell$, then
 $$\Tate_{\wp}(A) = \varprojlim_{r} A[\wp^r] = \Tate_\ell(A)\tensor_{\T} \T_{\wp}$$ is
    infinite. (For more details, see the proof of
\cite[Prop.~3.2]{ribet-stein}.)
\end{proof}

The analogues of Lemmas~\ref{lem:jtwofactors}--\ref{lem:somefactor}
below are true, with the same proofs, for $\m$ any Eisenstein prime.
We state and prove them for the $2$-Eisenstein prime, since that is
the main case of interest to us.  Let $\Jtwo$ be the {\em
$2$-Eisenstein quotient of~$J$}, where $\Jtwo$ is as defined in
\cite[\S{}II.10]{mazur}.  More precisely, we have the following:
\begin{lemma}\label{lem:jtwofactors}
  The simple factors of $\Jtwo$ correspond to the
  $\Gal(\Qbar/\Q)$-conjugacy classes of newforms~$f$ such that
  $A_f[\m]\neq 0$ (or equivalently, $A_f^{\vee}[\m]\neq 0$).
\end{lemma}
\begin{proof}
  On page 97 of \cite{mazur} we find that the $\C$-simple factors of
  $\Jtwo$ are in bijection with the irreducible
  components $\Spec(I_f)$ of $\Spec(\T)$ which meet the support of the
  ideal~$\m$, so the $I_f$ are the newform ideals contained in~$\m$. 
We have for any $I_f$,
$$
 I_f \subset \m \iff \T_m/(I_f)_\m \neq 0 \iff{}\Tate_\m(A_f)\neq 0
\iff{}A_f[\m]\neq 0.
$$
Note that the same argument applies to $A_f^{\vee}$.
\end{proof}

\begin{lemma}\label{lem:forcedtorsion}
Suppose $A$ and $B$ are abelian varieties equipped with an action
of the Hecke ring~$\T$ and that $\vphi:A\to B$ is a $\T$-module
isogeny.  If $\wp\subset \T$ is a maximal ideal and $B[\wp]\neq 0$,
then also $A[\wp]\neq 0$.
\end{lemma}
\begin{proof}
  Let $\psi:B\to A$ be the isogeny complementary to $\vphi$, so $\psi$
  is the unique isogeny such that $\psi\circ \vphi$ is multiplication
  by $\deg(\vphi)$.  Then $\psi$ is also a $\T$-module homomorphism
  (one can see this in various ways; one way is to use the rational
  representation on homology to view the endomorphisms as matrices
  acting on lattices, and to note that if matrices $M$ and $N$
  commute, then $M^{-1}$ and $N$ also commute).  By
  Lemma~\ref{lem:powkill}, the union $B[\wp^{\infty}]$ is infinite, so
  $\psi(B[\wp^{\infty}])\neq 0$.  Since $\psi(B[\wp^{\infty}])\subset
  A[\wp^{\infty}]$, Lemma~\ref{lem:powkill} implies that $A[\wp]\neq
  0$, as claimed.
\end{proof}

\begin{lemma}\label{lem:somefactor}
  Suppose $B\subset J_0(p)$ is a sum of abelian subvarieties $A_f$
  attached to newforms.  If $B[\m]\neq 0$, then there is some
  $A_f\subset B$ such that $A_f[\m]\neq 0$.
\end{lemma}
\begin{proof}
  There is something to be proved because if $x\in B[\m]$ it could be
  the case that $x=y+z$ with $y\in A_f$ and $z\in A_{g}$, but $x\not
  \in A_h$ for any~$h$.  Let $C = \oplus A_f$, where the $A_f\subset
  J_0(p)$ are simple abelian subvarieties of~$B$ corresponding to
  conjugacy classes of newforms.  Then there is an isogeny $\vphi:C\to
  B$ given by $$\vphi(x_1,\ldots,x_n) = x_1+\cdots +x_n,$$ where the
  sum is in $B\subset J_0(p)$.  By Lemma~\ref{lem:forcedtorsion},
  $C[\m]\neq 0$.  Since $C[\m]\isom \oplus A_f[\m]$, it follows that
  $A_f[\m]\neq 0$ for some~$A_f\subset B$.
\end{proof}

\begin{lemma}\label{lem:jtwobig}
If $4\mid n$, then $\dim \Jtwo > 1$.
\end{lemma}
\begin{proof}
  This follows from the remark on page 163 of \cite{mazur}.  Since the
  proof is only sketched there, we give further details for the
  convenience of the reader.  Because $4\mid n$, the cuspidal
  subgroup~$C$, which is generated in $J_0(p)$ by $(0)-(\infty)$ and
  is cyclic of order~$n$, contains an element of order~$4$.  Let
  $C(2)$ be the $2$-primary part of~$C$, and let $D =
  \ker(J_0(p)\to\Jtwo)$.  If there is a nonzero element in the kernel
  of the homomorphism $C(2) \to \Jtwo$, then $D[\m]\neq 0$, where $\m$
  is the $2$-Eisenstein prime.  But then by
  Lemma~\ref{lem:somefactor}, there is an $A_f\subset D$ such that
  $A_f[\m]\neq 0$.  By Lemma~\ref{lem:jtwofactors}, $A_f^{\vee}$ is a
  quotient of $\Jtwo$, so $A_f\subset (\Jtwo)^{\vee}$ so $A_f$ cannot
  be in $D$.  This contradiction shows that the map $C(2)\to \Jtwo$ is
  injective, so $\Jtwo$ contains a rational point of order~$4$.
  However, as mentioned in the introduction, $E_0(\Q)$ has order~$2$,
  so $\Jtwo\neq E_0$.  Thus $\Jtwo$ has dimension bigger than~$1$.
\end{proof}
Having established the above lemmas, we are now ready to deduce the
theorem.
\begin{proof}[Proof of Theorem~\ref{thm:main}]
It seems more straightforward to prove the equivalent
statement that the modular degree is even if
and only if $u\con 7\pmod{8}$, so we will prove
this instead.\vspace{1ex}\\
\noindent{$(\Longrightarrow)$ \em $u\con 7\pmod{8}$ implies that the
modular degree is even:} Writing $u=8k+7$ we see that $p =(8k+7)^2 +
64\equiv 1\pmod{16},$ so $16\mid (p-1)$ hence $4\mid n$.
By Lemma~\ref{lem:mtwo} and Lemma~\ref{lem:jtwofactors}, $E_0$ is a
factor of $\Jtwo$. By Lemma~\ref{lem:jtwobig}, the dimension of
$\Jtwo$ is bigger than~$1$, so by Lemma~\ref{lem:jtwofactors} there is
an $A_f$ distinct from $E_0$ such that $A_f[\m]\neq 0$.  Since
$A_f\subset B = \ker(J_0(p)\to E_0)$, it follows that $B[\m]\neq 0$.
As discussed on page 38 of \cite{mazur}, $J[\m]$ has dimension~$2$
over $\F_2$ so $E_0[\m]=J[\m]$, hence $B[\m] \subset E_0[\m]$.  It
follows that $2\mid \#(B\intersect E_0)$, so $E_0$ has even modular
degree.
\vspace{1ex}

\noindent{\em $(\Longleftarrow)$ Modular degree even implies that $u
\con 7\pmod{8}$:} Suppose that the modular degree~$m$ of $E_0$ is
even.  Letting $B=\Ker(J_0(p)\to E_0)$, we have $$E_0\cap{}B\isom
\ker(E_0\to J_0(p)\to E_0),$$ so $\Psi := E_0\cap B=E_0[m]$.
Lemma~\ref{lem:mtwo} and our assumption that $m$ is even imply that
$$
 E_0[\m]= E_0[2] \subset E_0[m] =\Psi,
 $$ so $\Psi[\m]\neq 0$. Since $\Psi[\m]\neq 0$, and $\Psi\subset B$,
 we have $B[\m]\neq 0$.  By Lemma~\ref{lem:somefactor}, there is some
 $A_f\subset B$ such that $A_f[\m]\neq 0$.  Then by
 Lemma~\ref{lem:jtwofactors} we see that $A_f$ is an isogeny factor of
 $\Jtwo$.  Thus $\Jtwo$ has dimension bigger than~$1$.  If
$u= 8k+3$, then $p=(8k+ 3)^2+64\equiv 9\pmod{16}$,
so that $2\mid\mid n$.  However, when $2\mid\mid n$,
 \cite[Prop.~III.7.5]{mazur} implies that $\Jtwo=E_0$, 
which is false, so $u \con 7\pmod{8}$.
\end{proof}

\begin{remark}
  Frank Calegari observed that Lemma~\ref{lem:jtwobig} and its
  converse also follow from conditions (i) and (v) of Th\'eor\`eme~3
  of \cite{merel:weil}.
\end{remark}

\section{The Stevens Conjecture for Neumann--Setzer Curves is True}
\label{sec:stevens}
Let~$E$ be an arbitrary elliptic curve over~$\Q$ of conductor~$N$.
Stevens conjectured in \cite{stevens} that the optimal quotient
of $X_1(N)$ in the isogeny class of~$E$ is the curve in the isogeny
class of~$E$ with minimal Faltings height.  In this section we
explain why this conjecture is true when~$N$ is prime.

Let $p=u^2+64$ be prime and $E_1$ and $E_0$ be as in
Section~\ref{sec:intro}.  In this section we verify that the curve
$E_1$ has smaller Faltings height than $E_0$, then show that $E_1$ is
$X_1(p)$-optimal.  The Stevens conjecture asserts that the
$X_1(p)$-optimal curve is the curve of minimal Faltings height in an
isogeny class, so our results verify the conjecture for
Neumann--Setzer curves.  In fact, the Stevens conjecture is true for all
isogeny classes of elliptic curves of prime conductor.  For if~$E$ is
an elliptic curve of prime conductor, then by \cite{setzer} there is
only one curve in the isogeny class of~$E$, unless~$E$ is a
Neumann--Setzer curve or the conductor of~$E$ is 11, 17, 19, or~37.
When the isogeny class of $E$ contains only one curve, that curve is
obviously both $X_1$-optimal and of minimal Faltings height.  The
conjecture is also well-known to be true for curves of conductor $11$,
$17$, $19$, or $37$ (see \cite{stevens}).  We note that Vatsal
\cite{vatsal} has recently extended results of Tang \cite{tang} that
make considerable progress toward the Stevens conjecture, but his work
is not applicable to Neumann--Setzer curves.

\begin{lemma}\label{lemma:analcomp}
The curve $E_1$ has smaller Faltings height than~$E_0$.
\end{lemma}
\begin{proof}
By \cite[Thm.~2.3, pg.~84]{stevens} it is enough to
exhibit an isogeny from
$E_1$ to $E_0$ whose extension to N\'eron models
is \'etale.
Let $\varphi$ be the isogeny $E_0\to E_1$ of degree $2$ whose kernel is
the subgroup generated by the point whose coordinates are
$(u/4,-u/8)$ in terms of the Weierstrass equation (\ref{eqn:e0})
for $E_0$, which is a global minimal model for $E_0$.
The kernel of~$\varphi$ does not extend
to an \'etale group scheme over~$\Z$,
since its special fiber at~$2$ is not \'etale (it has only
one $\overline{\F}_2$-point), so the morphism on N\'eron models
induced by $E_0\to E_1$ cannot be \'etale, since kernels of \'etale
morphisms are \'etale.
By \cite[Lemma~2.5]{stevens}
 the dual isogeny $E_1\to E_0$ extends
to an \'etale morphism of N\'eron models.
\end{proof}

\begin{proposition}
The curve $E_1$ is $X_1(p)$-optimal.
\end{proposition}

\begin{proof}
By \cite[\S5, Lem.~3]{mestre-oesterle}, $E_0$ is
an optimal quotient of $X_0(p)$, so we have an injection
$E_0\hra J_0(p)$.
As in \cite[pg.~100]{mazur}, let
$\Sigma$ be the kernel of the functorial map $J_0(p)\to J_1(p)$
induced by the cover $X_1(p)\to X_0(p)$.  By
\cite[Prop.~II.11.6]{mazur}, $\Sigma$ is the Cartier
dual of the constant subgroup scheme $U$, which
turns out to equal $J_0(p)(\Q)_{\tor}$.  Because
$\#(E_0 \intersect U) = 2$ and $E_0[2]$ is self dual,
we have $\#(E_0\intersect \Sigma)=2$.
Thus the image of $E_0$ in $J_1(p)$ is
the quotient of $E_0$ by the subgroup generated by
the rational point of order $2$ (note that the Cartier
dual of $\Z/2\Z$ is $\mu_2=\Z/2\Z$).  This quotient is $E_1$,
so $E_1\subset J_1(p)$, which implies that $E_1$ is
an optimal quotient of $X_1(p)$, as claimed.
\end{proof}

\begin{remark}
 The above proposition could also be proved in a slightly different manner.
 The Faltings height of an elliptic curve is $\sqrt{2\pi/\Omega}$ where
 $\Omega$ is the volume of the fundamental parallelogram associated to the
 curve. When the conductor is prime, we have by \cite{abbes-ullmo} that the
 Manin constants for $X_0(p)$ and $X_1(p)$ are~1; this says that for a
 $G$-optimal curve~$E$,
 the period lattice generated by~$G$ has covolume equal to $\Omega_E$.
 Since the lattice generated by $\Gamma_1(p)$ is contained in the lattice
 generated by $\Gamma_0(p)$ (and thus has larger covolume), the Faltings
 height of the $X_1(p)$-optimal curve must be less than or equal to that
 of the $X_0(p)$-optimal curve. So if these two curves differ, the $X_1(p)$-optimal
 curve must have smaller Faltings height.
\end{remark}

\begin{remark}
 On page 12 of \cite{mazur:aws}, there is
 a ``To be removed from the final draft'' comment that asks
 (in our notation) whether $E_0$ is $X_0(p)$-optimal when
 $p\equiv 1$ (mod~16). This is already answered by \cite{mestre-oesterle},
 whereas here we go further and show additionally
  that $E_1$ is $X_1(p)$-optimal.
\end{remark}

\section{Conjectures}\label{sec:conj}
\subsection{Refinement of Theorem~\ref{thm:main}}
The following conjectural refinement of Theorem~\ref{thm:main}
is supported by the experimental data of \cite{watkins}.  It
is unclear whether the method of proof of Theorem~\ref{thm:main}
can be extended to prove this conjecture.
\begin{conjecture}
If $u\con 7\pmod{8}$,
then $2$ exactly divides the modular degree of $E_0$ if and
only if $u \con 7\pmod{16}$.
\end{conjecture}

We can note that the pattern seems to end here; for curves with
$u\con 15\pmod{16}$ the data give no further information about
the 2-valuation of the modular degree. For instance, with $u=-17$
we have that $[1,1,1,-2,16]$ has modular degree $2^3\cdot 3$,
while with $u=175$ the curve $[1,1,1,-634,-6484]$ has modular degree
$2^2\cdot 3^3\cdot 5\cdot 23$. Similarly, we have that $u=-33$
gives the curve $[1,-1,1,-19,68]$ with modular degree $2^5\cdot 3$,
while $u=127$ gives the curve $[1,1,1,-332,-2594]$ of modular degree
$2^2\cdot 3^2\cdot 5\cdot 43$.

\subsection{The Parity of the Modular Degree}
According to Cremona's tables \cite{cremona}, of the $29755$ new
optimal elliptic curve quotients of $J_0(N)$ with $N<8000$, a mere
$89$ have odd modular degree, which is less than 0.3\%.  There are
$52878$ non Neumann--Setzer curves in the database of
\cite{brumer-mcguinness} with prime conductor $N\leq 10^7$; of these
curves $4592$, or 8\%, have odd modular degree (see \cite{watkins}).
One reason that curves tend to have even modular degree is that for
many curves the modular parametrization factors through an
Atkin-Lehner quotient.  Note that the method of \cite{watkins} used to
compute the modular degree is rigourous when the level is prime
because by \cite{abbes-ullmo} the Manin constant is~1 when the level
is odd and square-free.

If $f(x)=(8x+3)^2+64$, then it is a well-known conjecture
(see \cite{hlconj} and e.g., \cite[\S{}A1]{guy})
that there are infinitely many primes of the form $f(n)$ for some integer~$n$,
thus we make the following conjecture.

\begin{conjecture}
There are infinitely many elliptic curves over~$\Q$ with odd modular degree.
\end{conjecture}

Our data suggest the following conjecture:
\begin{conjecture}
If~$E$ is an optimal elliptic curve quotient of $J_0(p)$ with
$p\not\con 3 \pmod{8}$ and~$E$ is not a Neumann--Setzer curve then the
modular degree of~$E$ is even or $p=17$.
\end{conjecture}
There are $23442$ Brumer-McGuinness (see \cite{brumer-mcguinness})
curves of conductor $37\le p \leq 10^7$ with $p\con 3 \pmod{8}$, of
which $11815$ have even functional equation, of which $7322$ have
rank~$0$, and $4589$ have odd modular degree.  The significance of the
data concerning the rank is that the second author has conjectured
that $2^r$ divides the modular degree, where~$r$ is the rank.

\begin{remark}
Instead of asking about divisibility by~$2$, one could ask about
divisibility by~$p$.  The first author and Frank Calegari make a
conjecture about discriminants of Hecke algebras in
\cite{calegari-stein} that implies that
the modular degree of an elliptic curve of prime conductor~$p$ is not
divisible by~$p$.  This conjecture agrees with our data.
\end{remark}

\begin{table}
\caption{Frequency of a prime dividing $\Sha$}\label{tab:freq}
\begin{center}
\begin{tabular}{|c|c|c|c|c|c|}\hline
\quad{}restriction \quad{} & \quad{}number\quad{}&
\quad{} $p=3$\quad{} & \quad{} $p=5$\quad{} &
 \quad{}$p=7$\quad{} & \quad{}$p=11$\quad\mbox{}\\\hline
$u\equiv 1\pmod{8}$& 25559 & 33.2\% & 16.9\% & 9.2\% & 3.0\% \\
$u\equiv 3\pmod{8}$& 25557 & 39.7\% & 20.3\% & 14.3\% & 8.4\% \\
$u\equiv 5\pmod{8}$& 25584 & 36.2\% & 18.5\% & 11.5\% & 5.0\% \\
$u\equiv 7\pmod{8}$& 25612 & 34.3\% & 20.3\% & 14.3\% & 8.2\% \\\hline
$u\equiv 0\pmod{3}$& 34009 & 36.0\% & 18.7\% & 12.1\% & 6.0\% \\
$u\equiv 1\pmod{3}$& 34032 & 35.2\% & 18.6\% & 11.5\% & 5.6\% \\
$u\equiv 2\pmod{3}$& 34271 & 36.3\% & 19.7\% & 13.3\% & 6.9\% \\\hline
$u\equiv 0\pmod{5}$& 34208 & 33.1\% & 18.0\% & 11.4\% & 5.4\% \\
$u\equiv 2\pmod{5}$& 33879 & 37.1\% & 19.5\% & 12.8\% & 6.5\% \\
$u\equiv 3\pmod{5}$& 34225 & 37.3\% & 19.5\% & 12.7\% & 6.5\% \\\hline
total&102312&35.8\%&19.0\%&12.3\%&6.2\%\\
Delaunay&&36.1\%&20.7\%&14.5\%&9.2\%\\\hline
\end{tabular}
\end{center}
\end{table}

\subsection{Shafarevich--Tate Groups of Neumann--Setzer Curves}
We consider the distribution of $\Sha$ in the Neumann--Setzer family
(and note that similar phenomena occur in 
the related families listed in \cite{sw:ants5}).  We
look at $u$ with $u^2+64$ prime and less than $2\cdot 10^{12}$.  We
now take $u$ to be positive, which thus replaces the restriction
that~$u$ be 3~mod~4.  The heuristics of \cite{delaunay} would seem to
give us an idea of how often we expect a given prime to divide
$\Sha$. For instance, since Neumann--Setzer curves have rank~0, the
prime~3 should divide $\Sha$ about $36.1\%$ of the time. However,
Table \ref{tab:freq} gives a slightly different story with effects
seen that depend on the various congruential properties of~$u$.


\begin{thebibliography}{1}
\bibitem[AL96]{abbes-ullmo} A. Abbes, E. Ullmo, {\it \`A propos de la
conjecture de Manin pour les courbes elliptiques modulaires} (French).
Compositio Math. {\bf 103} (1996), no.~3, 269--286.
\bibitem[BM90]{brumer-mcguinness} A. Brumer, O. McGuinness, {\it The
behavior of the Mordell-Weil group of elliptic curves.}
Bull. Amer. Math. Soc. (N.S.) {\bf 23} (1990), no.~2, 375--382, data
available online at \verb+http://modular.fas.harvard.edu/~oisin+
\bibitem[CS03]{calegari-stein} F.~Calegari, W.~Stein, {\it Conjectures
About Discriminants of Hecke Algebras of Prime Level}, submitted
(2003).
\bibitem[Cre]{cremona} J. E. Cremona, {\it Elliptic Curves of
conductor $\le 17000$}, electronic tables available online at 
\verb+http://www.maths.nott.ac.uk/personal/jec/ftp/data+
\bibitem[Del01]{delaunay} C. Delaunay, {\it Heuristics on
Tate-Shafarevitch Groups of Elliptic Curves Defined over {\bf Q}.}
Experiment. Math. {\bf 10} (2001), no.~2, 191--196.
\bibitem[Del83]{deligne} P. Deligne, {\it Preuve des conjectures de
Tate et de Shafarevitch (d'apr\`es G. Faltings).} (French).  Seminaire
Bourbaki, Vol. 1983/84. Ast\'erisque No. 121-122, (1985), 25--41.
\bibitem[DI95]{diamond-im} F.~Diamond and J.~Im, \emph{Modular forms
and modular curves}, Seminar on {F}ermat's {L}ast {T}heorem,
Providence, RI, 1995, pp.~39--133.
\bibitem[Eme01]{emerton} M. Emerton, {\it Optimal Quotients of Modular
Jacobians,} preprint (2001).
\bibitem[Frey87]{frey} G. Frey, {\it Links between solutions of
$A-B=C$ and elliptic curves.}  In {\it Number Theory (Ulm, 1987),}
edited by H. P. Schlickewei and E. Wirsing, 31--62, Lecture Notes in
Mathematics, 1380, Springer-Verlag, New York, 1989.
\bibitem[Guy94]{guy} R. K. Guy, {\it Unsolved problems in number
theory,} Springer-Verlag, 1994.
\bibitem[HL22]{hlconj} G.\thinspace{}H. Hardy, J. E. Littlewood, {\it
Some Problems of 'Partitio Numerorum.' III.  On the Expression of a
Number as a Sum of Primes.} Acta. Math.  {\bf 44} (1922), 1--70.
\bibitem[LO91]{ling-oesterle} S. Ling, J. Oesterl\'e, {The Shimura
subgroup of $J_0(N)$.}  Ast\'erisque {\bf 196--197} (1991), 171--203.
\bibitem[Maz77]{mazur} B. Mazur, {\it Modular curves and the
Eisenstein ideal.}  Inst. Hautes \'Etudes Sci. Publ. Math. {\bf 47}
(1977), 33--186.
\bibitem[Maz98]{mazur:aws} B. Mazur, {\it Three Lectures about the
Arithmetic of Elliptic Curves.}  Rough, unedited, and preliminary
notes from lectures given at the 1998 Arizona Winter School. Available
at \verb+http://swc.math.arizona.edu/notes/files/98MazurLN.ps+
\bibitem[Mer96]{merel:weil} L.~Merel, \emph{L'accouplement de {W}eil
entre le sous-groupe de {S}himura et le sous-groupe cuspidal de
${J}\sb 0(p)$}, J. Reine Angew. Math. \textbf{477} (1996), 71--115.
\bibitem[MO89]{mestre-oesterle}
J.-F. Mestre, J. Oesterl\'e, {\it Courbes de Weil semi-stables de
discriminant une puissance $m$-i\`eme} (French).
J. Reine Angew. Math. {\bf 400} (1989), 173--184.
\bibitem[Mum70]{mumford}
D.~Mumford, \emph{Abelian varieties}, Published for the Tata Institute of
  Fundamental Research, Bombay, 1970, Tata Institute of Fundamental Research
  Studies in Mathematics, No. 5.
\bibitem[Neu71]{neumann} O. Neumann, {\it Elliptische Kurven mit
vorgeschriebenem Reduktionsverhalten. I, II} (German).
Math. Nachr. {\bf 49} (1971), 107--123, {\bf 56} (1973), 269--280.
\bibitem[Ogg74]{ogg} A. Ogg, {\it Hyperelliptic modular curves.}
Bull. Soc. Math. France {\bf 102} (1974), 449--462.
\bibitem[RS01]{ribet-stein} K.\thinspace{}A. Ribet and
W.\thinspace{}A. Stein, \emph{Lectures on {S}erre's conjectures},
Arithmetic algebraic geometry (Park City, UT, 1999), IAS/Park City
Math. Ser., vol.~9, Amer. Math. Soc., Providence, RI, 2001,
pp.~143--232.
\bibitem[Set75]{setzer} B. Setzer, {\it Elliptic Curves of prime
conductor.}  J. London Math. Soc. (2), {\bf 10} (1975), 367--378.
\bibitem[Ste89]{stevens} G. Stevens, {\it Stickelberger elements and
modular parametrizations of elliptic curves.} Invent. Math. {\bf 98}
(1989), no.~1, 75--106.
\bibitem[SW02]{sw:ants5} W.\thinspace{}A. Stein, M. Watkins, {\it A
Database of Elliptic Curves---First Report.}  In {\it Algorithmic
number theory} (Sydney 2002), 267--275, edited by C. Fieker and
D. Kohel, Lecture Notes in Comput. Sci., 2369, Springer, Berlin, 2002.
\bibitem[Tan97]{tang} S.-L. Tang, {\it Congruences between modular
forms, cyclic isogenies of modular elliptic curves, and integrality of
$p$-adic $L$-function.}  Trans. Amer. Math. Soc. {\bf 349} (1997),
no.~2, 837--856.
\bibitem[Vat03]{vatsal} V. Vatsal, {\it Multiplicative subgroups of
$J_0(N)$ and applications to elliptic curves,} preprint (2003).
\bibitem[Wat02]{watkins} M. Watkins, {\it Computing the modular degree
of an elliptic curve,} Experiment. Math. {\bf 11} (2002), no.~4,
487--502.
\bibitem[Wil95]{wiles} A.\thinspace{}J. Wiles, {\it Modular elliptic
curves and Fermat's last theorem.}  Ann. of Math. (2) {\bf 141}
(1995), no.~3, 443--551.
\end{thebibliography}
\end{document}